\documentclass[a4paper,11pt]{amsart}
\usepackage{amscd}
\usepackage{amssymb}
\usepackage{amsmath}
\usepackage{amsthm}
\usepackage{epsf}
\newtheorem{thm}{Theorem}[section]

\newtheorem{lemma}[thm]{Lemma}
\newtheorem{prop}[thm]{Proposition}

\newtheorem{defn}[thm]{Definition}

\newtheorem{obser}[thm]{Observation}
\newtheorem{rem}[thm]{Remark}
\newtheorem{conj}[thm]{Conjecture}

\def\min{\operatorname{min}}

\def\c1{\operatorname{c_1}}
\def\c2{\operatorname{c_2}}

\def\PP{{\mathbf P}}

\def\L{{\mathcal L}}

\def\Fq{{\mathbb F}_q}
\def\+{\oplus}                   % direct sum
\def\*{\otimes}                  % tensor product
\begin{document}

\title{Higher weights of Grassmann codes in terms of properties of Schubert unions}
\author{Sudhir R. Ghorpade, Trygve Johnsen, Arunkumar R. Patil and Harish K. Pillai}  

\address{Sudhir R. Ghorpade\\ Department of Mathematics\\ Indian Institute of Technology Bombay\\ Powai, Mumbai 400076, India}
\email{srg@math.iitb.ac.in}
\address{Trygve Johnsen\\ Department of Mathematics\\ University of Troms{\o}
\\ 9037 Troms{\o}, Norway}
\email{Trygve.Johnsen@uit.no}
\address{Arunkumar R. Patil\\ Shri Guru Gobind Shinghji Institute of Engineering \& Technology\\ Vishnupuri, Nanded 431606, India}
\email{arun.ittb@gmail.com}
\address{Harish K. Pillai\\ Department of  Electrical Engineering\\ Indian Institute of Technology Bombay\\ Powai, Mumbai 400076, India}
\email{hp@ee.iitb.ac.in}

%\address{Dept. of Mathematics\\ 
% University of Bergen\\ Johs. Brunsgt 12\\ 5008 Bergen, Norway}
%\email{johnsen@mi.uib.no, andreask@mi.uib.no}
\begin{abstract}
We describe the higher weights of the Grassmann codes $G(2,m)$ over finite fields $\Fq$
in terms of properties of Schubert unions, and in each case we determine the weight
as the minimum of two explicit polynomial expressions in $q$.
\end{abstract}
\maketitle

%\bigskip
%  \begin{abstract}
%       We give a simple result about some higher weights of Grassmann codes.
%We give some upper bounds for other higher weights, using Schubert cycles.
%We show that there are Schubert cycles that do not compute higher weights.
%At the ende we consider some natural questions.
%  \end{abstract}
%
%\vspace{.75cm}

% Section 1
\section{Introduction}
\label{intro}

Let $G(2,m)$ be the Grassmann variety of $2$-dimensional subspaces of a fixed $m$-dimensional vector space $V$ over the field ${\Fq}$ with $q$ elements. 
 By the standard Pl\"ucker 
coordinates $G(2,m)$ is embedded into
$\PP^{k-1}$ as a non-degenerate smooth subvariety, where $k=$$m
\choose  2$.
Let $C(2,m)$ be a code with $k \times n$ generator matrix $M$ where the $n$ columns of $M$ viewed as vectors in ${\Fq}^{k}$ are the coordinate representatives of the points of $G(2,m)$ under the Pl\"ucker embedding. Thus the columns of $M$ are only determined up to a non-zero multiplicative constant, and the ordering of the columns is arbitrary.
Nevertheless, the word length $n$, the dimension $k$, and the higher weigths
$d_1,\cdots,d_k$ are uniquely determined and independent of the choice of column order
and multiplicative constants. It is well known that the word length (the number of ${\Fq}-$rational points on $G(2,m)$ is 
\begin{equation} \label{number2}
n=\frac{(q^m-1)(q^{m-1}-1)}{(q^2-1)(q-1)},
\end{equation}
and that the dimension is $k=$$m \choose 2$ as suggested by the choice of notation.
In the thesis \cite{P} a formula for the higher weights $d_i$ was given in terms of
properties of a certain Young diagram. 

Moreover it is well known that for $i=1,\cdots,k$
   $$d_i = n-J_i,$$ where $J_i$ is the maximal number of ${\Fq}$-rational points from $G(2,m)$ that you can find on a codimension $i$ linear subspace of  the Pl\"ucker projective space
$\PP^{k-1}$. In \cite{HJR} one conjectured that
the higher weights were computed by the so-called Schubert unions, in the sense
that for any $i \in \{0,1,2,\cdots,k\}$ the maximal number of ${\Fq}$-rational points from $G(2,m)$ that you can find  
on a codimension $i$ linear subspace of  the Pl\"ucker space
$\PP^{k-1}$ can always be found on a codimension $i$ subspace that intersects $G(2,m)$ in such a Schubert union.

In this paper we show that the formula given in \cite{P} is identical with the one
predicted by the conjecture concerning Schubert unions, and that Schubert unions thus compute the higher weights. 
Moreover we may then utilize a procedure (Proposition 4.6 of \cite{HJR}, and Proposition 5.3 of \cite{HJR2}) for computing the
optimal Schubert unions, in the sense that they contain the maximal number of ${\Fq}$-rational points
among all Schubert unions spanning a linear subspace of Pl\"ucker space of the same dimension (at least for all large enough $q$). In view of the result in \cite{P}, it is clear that this procedure,
which was described without proof in \cite{HJR} (while a proof was given in \cite{HJR2}), enables us to  compute the higher weights of $C(2,m)$ (at least for large $q$). 
With permission from the two other authors of \cite{HJR2} we then also give a  
version of the proof here (with cosmetic changes only). For each $m$ and $i$ this reduces the computation of the higher weight $d_i$ of $C(2,m)$ over ${\Fq}$ to taking the minimum value of two explicit polynomial expressions in $q$.

\section{Basic Description of Schubert Unions}
\label{SUdef}

In this section we recall the basic facts from \cite{HJR} about Schubert unions,
necessary for pur purpose. We recall the well known definition of a Schubert
variety $\alpha =(a_1,...,a_l)$ in the Grassmann variety $G(l,m)$ over 
a field
$F$, and describe unions of such varieties.. 
Let $B=\{e_1,...,e_m\}$ be a basis of a $m$-dimensional 
vector space $V$ over $F$, and let $A_i=Span$$\{e_1,...,e_i\}$ in 
$V$, 
for $i=1,...,m.$  Then $A_{1}\subset A_{2}\subset \ldots \subset 
A_{m}=V$ form a complete flag of subspaces of $V$.   With respect to 
the basis $B$ there 
is the following canonical cell decomposition of $G(l,m)$.

For a given $l$-subspace $W$ of $V$
form an $(l \times m)$-matrix $M_{W}$ where the rows form a set of 
basis
vectors for $W$, each row expressed in terms of the basis $B$. 

We choose a basis 
%$<w_{1},\ldots,w_{l}>$ 
for $W$ such that the 
matrix $M_{W}$  have reduced lower left triangular form, i.e. 
the last nonzero entry in each row is $1$, each of these $1$'s are 
the 
only nonzero entries in their column, and each of these $1$'s lie in 
a column to the right of the trailing $1$ in the previous row. 
The trailing $1$ in row $i$ is then in column $a_{i}(W)$ where 
$$
a_{i}(W)={\rm min}\{j\in \{1, \dots , l\}  : \dim(W\cap A_{j})=i\}.
$$
%Obviously $$1\leq a_{1}(W)<a_{2}(W)<...<a_{l}(W)\leq m.$$
For $\alpha=(a_1, \dots ,a_l) \in {\bf Z}^l$ with $ 1 \le a_1 < a_2 < \dots < a_l \le    m$ 
we define the cell 
$$
C_{\alpha}=\{W\in G(l,m) : a_{i}(W)=a_{i}  \text{ for }  i=1, \dots ,l\}.
$$
The ordered $l$-tuples $\alpha$ belong to the grid
$$
I(l,m)=\{\beta=(b_1, \dots ,b_l)\in {\bf Z}^l | 1 \le b_1 < b_2 
< \dots < b_l \le    m\}.
$$
This grid is partially ordered by
$\alpha\leq \beta$ if $a_{i}\leq b_{i}$ for $i=1, \dots ,l$.

For each $\alpha\in I(l,m) $ the Schubert variety $S_{\alpha}$ is 
defined as:

$$
S_{\alpha}=\{W \in G(l,m) : \dim(W \cap A_{a_i}) \ge i \text{ for } 
i=1 \dots ,l\}=\bigcup_{\beta\leq\alpha}C_{\beta}.$$
This is the well-known cell decomposition of  $S_{\alpha}$.
%inherits a cell-decomposition from $G(l,m)$.

Next, we choose coordinates for the Pl\"ucker space 
$\PP^{k-1}=\PP(\wedge^l V)$, with
respect to the chosen basis $B$. 
Our choice of Pl\"ucker coordinates are the maximal minors 
of the matrix $M_{W}$. % (with alternating signs chosen in a standard way). 
These minors are indexed as
%by $G_{G(l,m)}$, so these Pl\"ucker coordinates are denoted by 
$\{X_{\alpha}(W)| \alpha\in I(l,m)\}$.
 
\begin{defn} For any $\alpha\in I(l,m)$, let 
$I_{\alpha}(l,m)=\{\beta\in I(l,m)|\beta\leq \alpha\}.$ 
%and $H_{\alpha}=G_{G(l,m)} - G_{\alpha}.$
\end{defn}

%We observe: 
For any $\alpha\in I(l,m)$, it is readily seen that 
$$
S_{\alpha}=\{W \in G(l,m) : X_{\beta}(W)=0 \text{ for  all } \beta\in 
I(l,m) \setminus I_{\alpha}(l,m)\}.
$$
Also, observe that 
%\leftline{
for $\alpha = (m-l+1,...,m-1,m)$, we get
$S_{\alpha}=G(l,m)$ and $I_{\alpha}(l,m)=I(l,m)$.
% and $H_{\alpha}=\emptyset$.
%}

 \begin{defn}
For a subset $E$ of $G(l,m)\subset \PP (\wedge^l V)$, let 
$\L(E)$ be the linear span of $E$ in
the projective  Pl\"ucker space $\PP (\wedge^l V),$
and let $L(E)$ the linear span of the
affine cone of $E$ %in $\wedge^l V$, i.e., 
in the affine cone of the projective Pl\"ucker space $\PP (\wedge^l V)$.
 \end{defn}
 
We will consider finite intersections and finite unions of such 
Schubert varieties
$S_{\alpha}$ with respect to our fixed flag.
Set $\alpha_i=(a_{i,1},a_{i,2}, \dots ,a_{i,l})$, for $i=1 \dots ,s$. 
It is clear that:
  $\cap_{i=1}^s S_{\alpha_i} = S_{\gamma}$,
where $\gamma=(g_1, \dots ,g_l)$, and $g_j = \min\{a_{1,j}, \dots ,a_{s,j}\}$ for $j=1, \dots,l.$ 
%is the minimum of the set $\{a_{1,j},a_{2,j},...,a_{s,j}\},$ for $j=1, \dots ,l.$  
Thus the intersection of a 
finite set of Schubert varieties $S_{\alpha}$ is again a Schubert variety.
In particular $\dim L(\cap S_{\alpha_i})$ is equal to
the cardinality of  $I_{\gamma}(l,m)$. 

\begin{defn} \label{uniongrid}
Given any $U \subset I(l,m)$,  the Schubert union $S_U$ and its associated Schubert union grid $I_U$ are defined by
$$
S_U=\bigcup_{\alpha \in U} S_{\alpha} \quad \text{ and } \quad 
I_U=\bigcup_{\alpha \in U}I_{\alpha}(l,m).
$$ 
%and set $H_U=G_{G(l,m)}-G_U$.
\end{defn}
We observe: $S_U \subset S_V$ if and only if $I_U \subset I_V$.
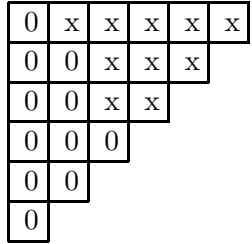
\begin{figure}[ht]
    \centering
     \setlength{\unitlength}{15pt}
\noindent
\begin{picture}(7,7)
  %   \thinlines
 %     \multiput(0,0)(1,0){7}{\line(0,1){7}}
 %  \multiput(0,0)(0,1){7}{\line(1,0){7}}
   \thicklines

   \put(-3,0){\line(0,1){6}}
  \put(-2,0){\line(0,1){6}}
 % \put(7,0){\line(0,1){7}}
%   \put(0,7){\line(1,0){6}}
   \put(-3,0){{\line(1,0){1}}}
  \put(-3,1){\line(1,0){2}}
 %  \put(0,0){{\line(1,0){1}}}
  \put(-1,1){\line(0,1){5}}
 % \put(1,1){{\line(1,0){1}}}
  \put(0,2){\line(0,1){4}}
  \put(-3,2){{\line(1,0){3}}}
  \put(1,3){\line(0,1){3}}
  \put(-3,3){{\line(1,0){4}}}
  \put(2,4){\line(0,1){2}}
  \put(-3,4){{\line(1,0){5}}}
  \put(3,5){\line(0,1){1}}
  \put(2,5){{\line(1,0){1}}}
 % \put(17,6){\line(0,1){1}}
  \put(-3,6){{\line(1,0){6}}}
   \put(0,3){{\line(0,1){2}}}
  \put(-3,5){{\line(1,0){6}}}
%  \put(4.4,4.3){0}  
\put(0.4,4.3){x}  
\put(1.4,4.3){x}  
  \put(-2.6,4.3){0}  
\put(0.4,5.3){x}  
\put(1.4,5.3){x}
\put(2.4,5.3){x}
%\put(-6.6,5.3){0}
\put(-2.6,5.3){0}
\put(-0.6,5.3){x}
\put(0.4,3.3){x}
\put(-1.6,5.3){x}
  \put(-1.6,4.3){0} 
 \put(-0.6,4.3){x}
  \put(-2.6,3.3){0}
  \put(-1.6,3.3){0} 
 \put(-0.6,3.3){x}
%  \put(1.1,0.0){Illustration of $I_U$, where $U=\{(1,7),(2,6),(3,4)\}$.}
 \put(-2.6,2.3){0}
  \put(-1.6,2.3){0} 
\put(-0.6,2.3){0}
 \put(-2.6,1.3){0}
  \put(-1.6,1.3){0}
\put(-2.6,0.3){0}
\end{picture}
   \caption{Illustration of $I_U$, where $U=\{(1,7),(2,6),(3,4)\}$.}
     \label{fig:u1}
\end{figure}

In the figure we have drawn the grid $I(2,7)=I_{(6,7)}(2,7)$ as the union of  $7 \choose  2$ squares arranged in a triangle. The bottom square is $(1,2)$, the top left square is $(1,7)$, and
the top right square is $(6,7)$.  The subset $I_U =I_{(1,7)}(2,7)\cup I_{(2,6)}(2,7)\cup I_{(3,4)}(2,7)$ of $I(2,7)$ is drawn as the union of the squares containing $0$'s.  For the special case $l=2$ it is clear that we can draw 
$I(2,m)=I_{(m-1,m)}(2,m)$ this way, 
with $m \choose 2$ squares arranged in a triangle,
for all $m \ge 2$. Then the bottom square is $(1,2)$, the top left square is $(1,m)$, and
the top right square is $(m-1,m)$.

We then make the following important observation:
\begin{obser} \label{divisionline}
For $l=2$, and any $m \ge 2$, draw any grid $I_U \subset I(2,m)$ as a union of squares with zeroes, and $I(2,m) - I_U$ as a union of squares with x's, as in the example. Then the division curve between the squares with $0$'s and the squares with x's is piecewise linear, and never goes upwards to the right.

Another formulation: The number $c_i$ of squares in $I_U$ from column number $i$
decreases strictly with $i$ until it reaches zero, and in each column of $I(2,m)$,
the squares included in $I_U$ are the $c_i$ ones closest to the bottom of that column.

Conversely: If we pick $c_1$ squares from column $1, \ldots, c_t$ squares from column
$t$, from the bottom for each column, for some $t \le m-1,$ and $c_1 > \cdots > c_t$,
then we obtain a subset of $I(2,m)$ which is equal to $I_U$ for some Schubert union.

\end{obser}

For more details, see \cite{HJR}, Section 3.
Since each Schubert variety $S_{\alpha}$ has a decomposition of cells 
$C_{\beta}$,
 and all finite intersections of these Schubert varieties are  again 
Schubert varieties, the union $S_{U}$ 
 also has a cell-decomposition inherited from $G(l,m)$:
$$S_{U}=\bigcup_{\alpha\in I_{U}}C_{\alpha}$$

The following is basically Proposition 2.3 from \cite{HJR}.

\begin{prop} \label{int-un}
Let $U \subset I(l,m)$ and let 
$S^U=\bigcap_{\alpha \in U} S_{\alpha},$ and $S_U=\bigcup_{\alpha \in U} S_{\alpha}.$
\begin{enumerate}
\item The intersection $S^U$
is itself a Schubert variety $S_{\gamma}$ with $I_{\gamma}(l,m)=\cap_{\alpha \in U} I_{\alpha}(l,m).$
\item  $\L(S_U)\cap G(l,m)=S_{U}.$
\item dim $L(S_{U})$ equals the cardinality of $I_{U}$.
\item  The number of ${\Fq}$-rational points on $S_{U}$ is 
 $g_U(q)=\Sigma_{(x_1,...x_l) \in I_{U}} q^{x_1+...+x_l-l(l+1)/2}$.
 \end{enumerate}
 \end{prop}
 We obtain  $g_U(q)=\Sigma_{(x_1,x_2) \in I_U}q^{x_1+x_2-3}$ for $l=2$.
It is well known that for the Grassmann codes $C(l,m)$ the higher weights $d_i$ satisfy $$d_i = n-J_i,$$ where $J_i$ is the the maximal number of ${\Fq}$-rational points from $G(l,m)$ on a codimension $i$ linear space in the Pl\"ucker space $\bf P$$^{k-1}$. Hence it gives meaning to say that $d_i$ is "computed by a Schubert union $U$" if a codimension $i$ space in Pl\"ucker space intersects $G(l,m)$ in a Schubert union where the number of ${\Fq}$-rational points is $J_i$. Our goal in the next section is to confirm Conjecture 5.4 of \cite{HJR} in the case $l=2$. This says:

\begin{conj} \label{Schubertmax} 
The higher weights of the Grassmann codes $G(l,m)$ are always computed 
by Schubert unions.
\end{conj}
 
 \section{Description of higher weights in terms of Young diagrams and Schubert unions }
\label{Young}

Let $Y_m$ be a set of boxes arranged in $m-1$ rows, with $m-i$ boxes in row number $i$ for $i=1,\cdots,m-1$. One inserts the number $2i+j-3$ in the $j$'th box from the left in the $i$'th row, for all $i,j$ in question. We display $Y_m$ in the following way
 (in the example $m=9$):

\begin{figure}[ht]
     \centering
     \setlength{\unitlength}{15pt}
\noindent
\begin{picture}(9,9)
  %   \thinlines
 %     \multiput(0,0)(1,0){7}{\line(0,1){7}}
 %  \multiput(0,0)(0,1){7}{\line(1,0){7}}
   \thicklines
   \put(0,0){\line(0,1){8}}
  \put(1,0){\line(0,1){8}}
 % \put(7,0){\line(0,1){7}}
%   \put(0,7){\line(1,0){6}}
   \put(0,0){{\line(1,0){1}}}
  \put(0,1){\line(1,0){2}}
   \put(0,0){{\line(1,0){1}}}
  \put(2,1){\line(0,1){7}}
 % \put(1,1){{\line(1,0){1}}}
  \put(3,2){\line(0,1){6}}
  \put(0,2){{\line(1,0){3}}}
  \put(4,3){\line(0,1){5}}
  \put(0,3){{\line(1,0){4}}}
  \put(5,4){\line(0,1){4}}
  \put(0,4){{\line(1,0){5}}}
  \put(6,5){\line(0,1){3}}
   \put(7,6){\line(0,1){2}}
    \put(8,7){\line(0,1){1}}
  \put(5,5){{\line(1,0){1}}}
 % \put(7,6){\line(0,1){1}}
  \put(0,6){{\line(1,0){7}}}
   \put(0,7){{\line(1,0){8}}}
      \put(0,8){{\line(1,0){8}}}
   \put(3,3){{\line(0,1){2}}}
  \put(0,5){{\line(1,0){6}}} 
\put(3.4,7.3){3}  
\put(4.4,7.3){4}
\put(5.4,7.3){5}
\put(0.4,7.3){0}
\put(2.4,7.3){2}
\put(1.4,7.3){1}
\put(6.4,7.3){6}
\put(7.4,7.3){7}
\put(3.4,6.3){5}  
\put(4.4,6.3){6}
\put(5.4,6.3){7}
\put(0.4,6.3){2}
\put(2.4,6.3){4}
\put(1.4,6.3){3}
\put(6.4,6.3){8}
\put(3.4,4.3){9}  
\put(4.1,4.3){10}  
  \put(0.4,4.3){6}  
\put(3.4,5.3){7}  
\put(4.4,5.3){8}
\put(5.4,5.3){9}
\put(0.4,5.3){4}
\put(2.4,5.3){6}
\put(3.1,3.3){11}
\put(1.4,5.3){5}
  \put(1.4,4.3){7} 
 \put(2.4,4.3){8}
  \put(0.4,3.3){8}
  \put(1.4,3.3){9} 
 \put(2.1,3.3){10}
%  \put(4.1,0.0){$Y_9$}
 \put(0.1,2.3){10}
  \put(1.1,2.3){11} 
\put(2.1,2.3){12}
 \put(0.1,1.3){12}
  \put(1.1,1.3){13}
\put(0.1,0.3){14}
\end{picture}
   \caption{The Young tableau $Y_9$}
     \label{fig:u2}
\end{figure}

A strict subtableau of $Y_m$ is an arrangement of a subset of boxes, where one picks $\lambda_i$ 
consecutive boxes from the left end of the $i$'th row, for $i=1,\cdots,t$, where
$t \le m-1$ and $\lambda_1 > \lambda_2> \cdots > \lambda_t$, and keep the numbers in each chosen box. One displays the subtableaux by letting all the chosen boxes stay where they are, and removing all the other boxes. The total number of boxes in the subtableau, viz., $\lambda_1+ \cdots + \lambda_t$ is called the \emph{area} of that subtableau. 
An example with $m=9$ is shown in the left half of the next figure.

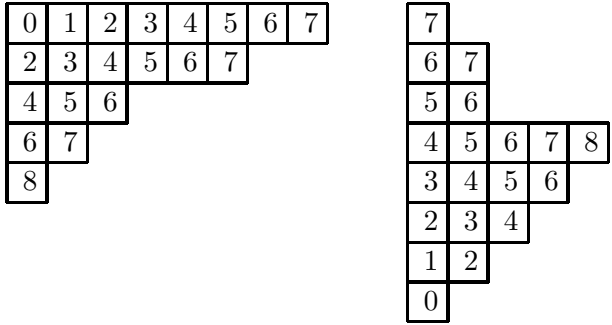
\begin{figure}[ht]
     \centering
     \setlength{\unitlength}{15pt}
\noindent
\begin{picture}(22,9)
  %   \thinlines
 %     \multiput(0,0)(1,0){7}{\line(0,1){7}}
 %  \multiput(0,0)(0,1){7}{\line(1,0){7}}
   \thicklines
   \put(0,3){\line(0,1){5}}
  \put(1,3){\line(0,1){5}}
   
  % \put(0,0){{\line(1,0){1}}}
  \put(2,4){\line(0,1){4}}
 % \put(1,1){{\line(1,0){1}}}
  \put(3,5){\line(0,1){3}}
 
  \put(4,6){\line(0,1){2}}
  \put(0,3){{\line(1,0){1}}}
  \put(5,6){\line(0,1){2}}
  \put(0,4){{\line(1,0){2}}}
  \put(6,6){\line(0,1){2}}
   \put(7,7){\line(0,1){1}}
    \put(8,7){\line(0,1){1}}
  %\put(5,5){{\line(1,0){1}}}
 % \put(7,6){\line(0,1){1}}
  \put(0,6){{\line(1,0){6}}}
   \put(0,7){{\line(1,0){8}}}
      \put(0,8){{\line(1,0){8}}}
  % \put(3,3){{\line(0,1){2}}}
  \put(0,5){{\line(1,0){3}}} 
\put(3.4,7.3){3}  
\put(4.4,7.3){4}
\put(5.4,7.3){5}
\put(0.4,7.3){0}
\put(2.4,7.3){2}
\put(1.4,7.3){1}
\put(6.4,7.3){6}
\put(7.4,7.3){7}
\put(3.4,6.3){5}  
\put(4.4,6.3){6}
\put(5.4,6.3){7}
\put(0.4,6.3){2}
\put(2.4,6.3){4}
\put(1.4,6.3){3}
%\put(6.4,6.3){8}
%\put(3.4,4.3){9}  
%\put(4.1,4.3){10}  
  \put(0.4,4.3){6}  
%\put(3.4,5.3){7}  
%\put(4.4,5.3){8}
%\put(5.4,5.3){9}
\put(0.4,5.3){4}
\put(2.4,5.3){6}
\put(1.4,5.3){5}
  \put(1.4,4.3){7} 
 %\put(2.4,4.3){8}
  \put(0.4,3.3){8}
%
%  \put(-3.1,1.0){$\lambda_1=8, \lambda_2=6, \lambda_3=3,\lambda_4=2, \lambda_5=1$}
%
   \put(10,0){\line(0,1){8}}
   \put(10,8){\line(1,0){1}}
    \put(10,7){\line(1,0){2}}
  \put(11,0){\line(0,1){8}}
 % \put(7,0){\line(0,1){7}}
%   \put(0,7){\line(1,0){6}}
   \put(10,0){{\line(1,0){1}}}
  \put(10,1){\line(1,0){2}}
   \put(10,0){{\line(1,0){1}}}
  \put(12,1){\line(0,1){6}}
 % \put(1,1){{\line(1,0){1}}}
  \put(13,2){\line(0,1){3}}
  \put(10,2){{\line(1,0){3}}}
  \put(14,3){\line(0,1){2}}
  \put(10,3){{\line(1,0){4}}}
  \put(15,4){\line(0,1){1}}
  \put(10,4){{\line(1,0){5}}}
  %\put(16,5){\line(0,1){1}}
  %\put(15,5){{\line(1,0){1}}}
 % \put(7,6){\line(0,1){1}}
  \put(10,6){{\line(1,0){2}}}
   \put(13,3){{\line(0,1){2}}}
  \put(10,5){{\line(1,0){5}}}
%\put(0.4,4.3){x}  
\put(13.4,4.3){7}  
\put(14.4,4.3){8}  
\put(10.4,4.3){4}  
%\put(13.4,5.3){x}  
%\put(14.4,5.3){x}
%\put(15.4,5.3){x}
%\put(0.4,5.3){x}
\put(10.4,5.3){5}
\put(10.4,6.3){6}
\put(10.4,7.3){7}
%\put(12.4,5.3){x}
\put(13.4,3.3){6}
\put(11.4,5.3){6}
\put(11.4,6.3){7}
  \put(11.4,4.3){5} 
 \put(12.4,4.3){6}
  \put(10.4,3.3){3}
  \put(11.4,3.3){4} 
 \put(12.4,3.3){5}
%  \put(14.1,0.0){$I_{(1,9)}(2,9)$$\cup I_{(2,8)}(2,9)$$\cup I_{(5,6)}(2,9)$ }
%is the dual of $S_{(3,5)}$}
 \put(10.4,2.3){2}
  \put(11.4,2.3){3} 
\put(12.4,2.3){4}
 \put(10.4,1.3){1}
  \put(11.4,1.3){2}
\put(10.4,0.3){0}
\end{picture}
 \caption{A subtableau of $Y_9$ with $\lambda_1=8, \lambda_2=6, \lambda_3=3,\lambda_4=2, \lambda_5=1$ and the corresponding Schubert union grid
 $I_{(1,9)}(2,9) \cup I_{(2,8)}(2,9) \cup I_{(5,6)}(2,9)$}
     \label{fig:u3}
\end{figure}

In \cite{P}(Theorem 6.18) one gives the following result.
\begin{thm} \label{Yg}
Fix an integer $r$, where $1 \le r\le k$. Let ${\mathcal F}=\{T_1,T_2,\cdots,T_p\}$
be the family of distinct strict subtableaux with area $r$ of $Y$. Let $a_{d,h}$ denote the number of times that the number $d$ occurs in the $h$th subtableaux $T_h \in {\mathcal F}.$ Associate $\gamma_h=\Sigma_{i}a_{d,h}q^i$ to the subtableaux $T_h$ for each $h$. Then for every $r$, we have $d_{k-r}=n-max_h\gamma_h.$
\end{thm}

The simple observation is now: There is a bijection between the set of boxes in $Y_m$
and the set of boxes in the grid $I(2,m)$ as follows: The $j$'th box from the left in the $i$'th row of $Y$
corresponds to the $j$'th box from the bottom in the $i$'th column of $I(2,m)$, for $i=1,\cdots,m-1$, and all $j$ in question. Then 
the number inserted in any chosen box $B$ of $Y$ is equal to the Krull dimension of
the Schubert variety $S_{\alpha}$ for the box $\alpha_B$ in $I(2,m)$ that
corresponds to $B$ under this bijection. In the right half of the figure above  the grid $I_U$
associated to the Schubert union $S_{(1,9)}$$\cup S_{(2,8)}$$\cup S_{(5,6)}$
corresponding to $(\lambda_1,\lambda_2,\lambda_3,\lambda_4,\lambda_5)=(8,6,3,2,1)$
for $m=9$, is displayed. The Krull dimensions of the $S_{\alpha}$ for each box ${\alpha}$ are included.

This immediately gives the crucial:
\begin{lemma} \label{basic}
There is a one-to-one correspondence between the set of strict subtableaux of $Y$ and
the set of subgrids $I_U$ of $I(2,m)$ for all possible Schubert unions $S_U$ (with respect to a fixed chosen flag). Under this correspondence the number $\gamma_h$
associated to a subtableaux $T_h$ is equal to the number $g_U$ of ${\Fq}$-rational points
in the Schubert union $S_U$.
\end{lemma} 

The results of Theorem \ref{Yg} and Lemma \ref{basic}  taken together confirm Conjecture \ref{Schubertmax} in the case $l=2$, and we obtain the main result of this paper:

\begin{thm} \label{Schubertmaxthm} 
The higher weights of the Grassmann codes $G(2,m)$ are always computed 
by Schubert unions.
\end{thm}

\section{Schubert unions with a maximal number of points}
\label{algo}

In this section we will show for each $m \ge 2$ how to find the Schubert unions with the largest
number of ${\Fq}$-rational points among those Schubert unions of fixed spanning dimension
$K=\dim L(S_U)$, for each $K=0,\dots,k$.

\begin{defn} \label{orders} 
{\rm Fix a dimension $0\leq K\leq {m\choose 2}$, and
consider the set of Schubert unions $\{S_U\}_{K}$ in $G(l,m)$ with 
spanning dimension $K$.

Then we order the elements in $\{S_U\}_{K}$ according to the 
lexicographic 
order on the polynomials $g_U$. In other words $S_U > S_V$ if $\deg g_U > 
\deg g_V$ or $\deg g_U = \deg g_V$, and the coefficient of $g_U$ is 
larger than that of $g_V$ in the largest degree where the 
coefficients differ. 
We call this the order with respect to $g_U$. }
\end{defn}

Furthermore we need:

\begin{defn} \label{thresholds} 
{\rm Set $\nu_i=1+\cdots+i$, and $\mu_i=(m-1)+\cdots+(m-i)$ for $i=1,\ldots,m-1$ and set $\nu_0=\mu_0=0.$
Now for $1\le K \le$$ m \choose  2$ there exists a unique $x$ such that $\mu_x<K \le\mu_{x+1}$, and
a unique $z$ such that $\nu_z<K \le\nu_{z+1}$ ; we set $S_L=S_{(x,m)} \bigcup S_{(x+1,K-\mu_{x}+x+1)}$
and $S_R=S_{(z,z+1)} \bigcup S_{(K-\nu_z,z+2)}.$ }
%(with the convention $S_{(0,b)}=S_{(a,0)}=\emptyset$ for all $a,b$.)}
\end{defn}
\begin{prop} \label{leftright} 
 Fix a dimension $0\leq K\leq {m\choose 2}.$ 
Then $S_{L}$ or $S_{R}$ is maximal in  $\{S_U\}_{K}$  with respect 
to the natural lexicographic order on the polynomials $g_U$. 
Furthermore, the one(s) that is(are) maximal with respect to $g_U$, 
also has(have) the maximum number of points over ${\Fq}$ for all large 
enough $q$.  
\end{prop}

\begin{rem} \label{both}
{\rm It is clear that for given $m$ and $K$, we obtain $I_U$
of the described $S_L$ by ``filling up as many 
columns of the $I(2,m)$-grid as we can from the left and filling in the remaining 
boxes (if any) from the bottom in the next column''. 

Likewise we obtain the $I_U$ of the described $S_R$ by 
``filling up as many rows of the $I(2,m)$-grid as we can from 
the bottom and filling in the remaining boxes (if any) from the left in the next row''.}
\end{rem}

\begin{proof}
 Given the spanning dimension $K$, let $d=d(K)$ be the
maximal Krull dimension for the Schubert unions $\{S_U\}_{K}$. 
This Krull dimension is the crucial ingredient in 
our argument, since the Krull dimension is the degree of the 
polynomial $g_{U}$. We will find the maximal polynomial $g_{U}$ 
in the lexicographic order.
The fact that the union(s) that is(are) maximal with respect to 
$g_U$, also has(have) the maximum number of points over ${\Fq}$ for all 
large enough $q$, is obvious.
Our argument is visualized by $I(2,m)$, arranged
as a set of squares in a triangle as on the next figure.
%defined by $$G_{(m-1,m)}=\{(x,y)| 1\leq x<y\leq m\}.$$
Each point $(a,b)\in I(2,m)$ defines a Schubert variety $S_{(a,b)}$ 
with Krull dimension $d(a,b)=a+b-3$.  Therefore the Schubert varieties 
with a fixed Krull dimension lie on the diagonal 
$$D_{d}=\{(x,y) |\quad 1\leq x<y\leq m, \quad x+y-3=d\}.$$
Let as above $$I_{(a,b)}(2,m)=\{(x,y)\in I(2,m)| x\leq a, y\leq b\},$$
and $$I_{U}=\bigcup_{(a,b)\in U}I_{(a,b)}(2,m).$$

By definition of $d=d(K)$, there is a
Schubert union $S_U$ of spanning dimension $K$ with an $I_U$ that 
contains 
a point $(a,b)$ on the diagonal $D_{d}$, i.e. $a+b-3=d$, but there is 
no such union $S_U$
with $I_U$-grid that contains a point on the diagonal  $D_{d+1}$. 

The {\it cardinality} $c(x,y)$ of a $I_{(x,y)}(2,m)$ defines the 
function
$$c:I(2,m)\to Z, \quad (x,y)\mapsto xy- \frac{x(x+1)}{2}.$$

The restriction of this function to the diagonal $D_{d}$ is defined by
$$c(x,d-x+3)=x(d-x+3)-\frac{x(x+1)}{2}, \quad {\rm for}\quad {\rm  
max} \{d+2-m,0\} < x < \frac{d+3}{2}$$
which is clearly quadratic and concave.  Therefore it attains its 
minimum $C(d)$, when $x$ is minimal or maximal, i.e. at one of the 
end points of 
the diagonal $D_{d}$.

  Clearly $$C(d(K)) \le K \le C(d(K)+1)-1.$$ 

 We say that a point $(a,b)\in D_{d(K)}$ is {\bf admissible}, if 
                $$c(a,b)< C(d+1),$$
 i.e. has less cardinality than any point in the next diagonal. Equivalently, 
  $(a,b)\in D_{d}$ is admissible if $I_{(a,b)}(2,m)\subset I_{U}$ for some Schubert union $S_U$ of spanning dimension at least $K$. For us the crucial fact is that for $(a,b)\in D_{d}$, being admissible is a necessary condition for having $I_{(a,b)}(2,m)\subset I_{U}$ for some Schubert union $S_U$ of spanning dimension exactly $K$. 
 
 Next, we characterize the admissible points by which diagonal 
$D_{d}$ 
 they belong to.

\begin{lemma} \label{admiss}
    Consider the diagonal 
     $$D_{d}=\{(x,y)|x+y-3=d\}=\{(x,d-x+3)\quad |\quad {\rm 
max}\{d+2-m,0\}<x<\frac{d+3}{2}\}$$
    (i) Let $d\leq m-3$, then the only admissible point on $D_{d}$ is 
    $(1,d+2)$, except when $d=2$, where $(2,3)$ is also 
    admissible.  
   
(ii) Let $d>m-3$, then $(x,d-x+3)$ is an admissible point on the 
diagonal $D_{d}$
only if $d+3-m\leq x\leq d+4-m$ or $\frac{d}{2}\leq 
x\leq\frac{d+2}{2}$, 
%i.e. only if it is among the two points with the smallest value of $x$, or the 
%two points with largest value $x$, 
with one exception, namely when $m=11$ and $d=10$, then the point $(4,9)$ is also admissible.

(iii) If  $m>10$, then  the point $(x,m)$ is admissible, only if $x 
\ge m-3$ or $x \le
\frac{m}{5} + 2$ if $x+m$ is odd, and only if $x \ge m-3$  or  $x \le
\frac{m}{5} + 1$ if $x+m$ is even.

(iv) If  $m>10$, then  the point $(x,m-1)$ is admissible only if $x \ge 
m-4$ or
$x \le \frac{m}{5} + 1$ if $x+m$ is odd, and only if $x \ge m-4$ or 
$x \le
\frac{m}{5} + 2$ if $x+m$ is even.
\end{lemma}
\begin{rem} \label{interior}
{\rm 
(a) The lemma implies, in very rough terms, that apart from the
leftmost column, the two
uppermost rows, and the two right-lowest points on each diagonal
(it's really only necessary to consider the right half of them, in
addition to $(2,3)$) then the ``interior''  that
remains contains no admissable points (except $(4,9)$ for $m=11$).

(b) If $d$ is odd, then part (ii) of the lemma implies that the point 
of the diagonal $D_{d}$ with the next to largest $x$, is 
non-admissible, so that the admissible ones must be among the two to 
the left, and the one to the right.}
\end{rem}

As an illustration we now show the grid $I(2,15)$ with the cost values
included for each element in the grid:

\begin{figure}[ht]
     \centering
     \setlength{\unitlength}{15pt}
\noindent
\begin{picture}(15,15)
  %   \thinlines
 %     \multiput(0,0)(1,0){7}{\line(0,1){7}}
 %  \multiput(0,0)(0,1){7}{\line(1,0){7}}
% \put(4.1,0.0){$G_{(2,15)}$}
    \thicklines
   \put(0,0){\line(0,1){14}}
  \put(1,0){\line(0,1){14}}
   \put(0,0){{\line(1,0){1}}}
  \put(0,1){\line(1,0){2}}
   \put(0,0){{\line(1,0){1}}}
  \put(2,1){\line(0,1){13}}
 % \put(1,1){{\line(1,0){1}}}
  \put(3,2){\line(0,1){12}}
  \put(0,2){{\line(1,0){3}}}
  \put(4,3){\line(0,1){11}}
  \put(0,3){{\line(1,0){4}}}
  \put(5,4){\line(0,1){10}}
  \put(0,4){{\line(1,0){5}}}
  \put(6,5){\line(0,1){9}}
  \put(5,5){{\line(1,0){1}}}
 % \put(7,6){\line(0,1){1}}
  \put(0,6){{\line(1,0){7}}}
   \put(3,3){{\line(0,1){10}}}
  \put(0,5){{\line(1,0){6}}} 
\put(0,7){{\line(1,0){8}}}
\put(5,6){{\line(0,1){8}}}
\put(6,7){{\line(0,1){7}}}
\put(7,6){{\line(0,1){8}}}
\put(8,7){{\line(0,1){7}}}
\put(9,8){{\line(0,1){6}}}
\put(10,9){{\line(0,1){5}}}
\put(11,10){{\line(0,1){4}}}
\put(12,11){{\line(0,1){3}}}
\put(13,12){{\line(0,1){2}}}
\put(14,13){{\line(0,1){1}}}
\put(0,14){{\line(1,0){14}}}
\put(0,8){{\line(1,0){9}}}
\put(0,9){{\line(1,0){10}}}
\put(0,10){{\line(1,0){11}}}
\put(0,11){{\line(1,0){12}}}
\put(0,12){{\line(1,0){13}}}
\put(0,13){{\line(1,0){14}}}
\put(3.2,4.3){14}  
\put(4.2,4.3){15}  
  \put(0.4,4.3){5}  
\put(3.1,5.3){18}  
\put(4.1,5.3){20}
\put(5.1,5.3){21}
\put(0.4,5.3){6}
\put(2.1,5.3){15}
%\put(2,5){{\line(1,1){1}}}
\put(3.1,3.3){10}
\put(1.1,5.3){11}
  \put(1.4,4.3){9} 
 \put(2.1,4.3){12}
  \put(0.4,3.3){4}
  \put(1.4,3.3){7} 
 \put(2.4,3.3){9}
%  \put(4.1,0.0){$I(2,15)$}
 \put(0.4,2.3){3}
  \put(1.4,2.3){5} 
\put(2.4,2.3){6}
 \put(0.4,1.3){2}
  \put(1.4,1.3){3}
\put(0.4,0.3){1}
\put(0.4,6.3){7}
\put(0.4,7.3){8}
\put(0.4,8.3){9}
\put(0.1,9.3){10}
\put(0.1,10.3){11}
\put(0.1,11.3){12}
\put(0.1,12.3){13}
\put(0.1,13.3){14}
\put(1.1,6.3){13}
\put(1.1,7.3){15}
\put(1.1,8.3){17}
\put(1.1,9.3){19}
\put(1.1,10.3){21}
\put(1.1,11.3){23}
\put(1.1,12.3){25}
\put(1.1,13.3){27}
\put(2.1,6.3){18}
\put(2.1,7.3){21}
\put(2.1,8.3){24}
\put(2.1,9.3){27}
\put(2.1,10.3){30}
\put(2.1,11.3){33}
\put(2.1,12.3){36}
\put(2.1,13.3){39}
\put(3.1,6.3){22}
\put(3.1,7.3){26}
\put(3.1,8.3){30}
\put(3.1,9.3){34}
\put(3.1,10.3){38}
\put(3.1,11.3){42}
\put(3.1,12.3){46}
\put(3.1,13.3){50}
\put(4.1,6.3){25}
\put(4.1,7.3){30}
\put(4.1,8.3){35}
\put(4.1,9.3){40}
\put(4.1,10.3){45}
\put(4.1,11.3){50}
\put(4.1,12.3){55}
\put(4.1,13.3){60}
\put(5.1,6.3){27}
\put(5.1,7.3){33}
\put(5.1,8.3){39}
\put(5.1,9.3){45}
\put(5.1,10.3){51}
\put(5.1,11.3){57}
\put(5.1,12.3){63}
\put(5.1,13.3){69}
\put(6.1,6.3){28}
\put(6.1,7.3){35}
\put(6.1,8.3){42}
\put(6.1,9.3){49}
\put(6.1,10.3){56}
\put(6.1,11.3){63}
\put(6.1,12.3){70}
\put(6.1,13.3){77}
\put(7.1,7.3){36}
\put(7.1,8.3){44}
\put(7.1,9.3){52}
\put(7.1,10.3){60}
\put(7.1,11.3){68}
\put(7.1,12.3){76}
\put(7.1,13.3){84}
\put(8.1,8.3){45}
\put(8.1,9.3){54}
\put(8.1,10.3){63}
\put(8.1,11.3){72}
\put(8.1,12.3){81}
\put(8.1,13.3){90}
\put(9.1,9.3){55}
\put(9.1,10.3){65}
\put(9.1,11.3){75}
\put(9.1,12.3){85}
\put(9.1,13.3){95}
\put(10.1,10.3){66}
\put(10.1,11.3){77}
\put(10.1,12.3){88}
\put(10.1,13.3){99}
\put(11.1,11.3){78}
\put(11.1,12.3){90}
\put(11.0,13.3){\small{102}}
\put(12.1,12.3){91}
\put(12.0,13.3){\small{104}}
\put(13.0,13.3){\small{105}}
\end{picture}
   \caption{The grid $I(2,15)$ with the cost values
included for each element}
     \label{fig:u4}
\end{figure}

As an additional illustration we show the same grid with the
non-admissible points crossed out:
\begin{figure}[ht]
     \centering
     \setlength{\unitlength}{15pt}
\noindent
\begin{picture}(15,15)
  %   \thinlines
 %     \multiput(0,0)(1,0){7}{\line(0,1){7}}
 %  \multiput(0,0)(0,1){7}{\line(1,0){7}}
   \thicklines
%   \put(4.1,0.0){$I(2,15)$ with non-admissible points crossed out}
    \put(0,0){\line(0,1){14}}
  \put(1,0){\line(0,1){14}}
   \put(0,0){{\line(1,0){1}}}
  \put(0,1){\line(1,0){2}}
   \put(0,0){{\line(1,0){1}}}
  \put(2,1){\line(0,1){13}}
 % \put(1,1){{\line(1,0){1}}}
  \put(3,2){\line(0,1){12}}
  \put(0,2){{\line(1,0){3}}}
  \put(4,3){\line(0,1){11}}
  \put(0,3){{\line(1,0){4}}}
  \put(5,4){\line(0,1){10}}
  \put(0,4){{\line(1,0){5}}}
  \put(6,5){\line(0,1){9}}
  \put(5,5){{\line(1,0){1}}}
 % \put(7,6){\line(0,1){1}}
  \put(0,6){{\line(1,0){7}}}
   \put(3,3){{\line(0,1){10}}}
  \put(0,5){{\line(1,0){6}}} 
\put(0,7){{\line(1,0){8}}}
\put(5,6){{\line(0,1){8}}}
\put(6,7){{\line(0,1){7}}}
\put(7,6){{\line(0,1){8}}}
\put(8,7){{\line(0,1){7}}}
\put(9,8){{\line(0,1){6}}}
\put(10,9){{\line(0,1){5}}}
\put(11,10){{\line(0,1){4}}}
\put(12,11){{\line(0,1){3}}}
\put(13,12){{\line(0,1){2}}}
\put(14,13){{\line(0,1){1}}}
\put(0,14){{\line(1,0){14}}}
\put(0,8){{\line(1,0){9}}}
\put(0,9){{\line(1,0){10}}}
\put(0,10){{\line(1,0){11}}}
\put(0,11){{\line(1,0){12}}}
\put(0,12){{\line(1,0){13}}}
\put(0,13){{\line(1,0){14}}}
\put(3.2,4.3){14}  
\put(4.2,4.3){15}  
  \put(0.4,4.3){5}  
\put(3.1,5.3){18}  
\put(4.1,5.3){20}
\put(5.1,5.3){21}
\put(0.4,5.3){6}
\put(2.1,5.3){15}
\put(3.1,3.3){10}
\put(1.1,5.3){11}
  \put(1.4,4.3){9} 
 \put(2.1,4.3){12}
  \put(0.4,3.3){4}
  \put(1.4,3.3){7} 
 \put(2.4,3.3){9}
 \put(0.4,2.3){3}
  \put(1.4,2.3){5} 
\put(2.4,2.3){6}
 \put(0.4,1.3){2}
  \put(1.4,1.3){3}
\put(0.4,0.3){1}
%\put(0.4,-1.7){$d$}
{\small 
\put(1.4,-0.7){0}
\put(1.9,-0.2){1}
\put(2.4,0.3){2}
\put(2.9,0.8){3}
\put(3.4,1.3){4}
\put(3.9,1.8){5}
\put(4.4,2.3){6}
\put(4.9,2.8){7}
\put(5.4,3.3){8}
\put(5.9,3.8){9}
\put(6.4,4.3){10}
\put(6.9,4.8){11}
\put(7.4,5.3){12}
\put(7.9,5.8){13}
\put(8.4,6.3){14}
\put(8.9,6.8){15}
\put(9.4,7.3){16}
\put(9.9,7.8){17}
\put(10.4,8.3){18}
\put(10.9,8.8){19}
\put(11.4,9.3){20}
\put(11.9,9.8){21}
\put(12.4,10.3){22}
\put(12.9,10.8){23}
\put(13.4,11.3){24}
\put(13.9,11.8){25}
\put(14.4,12.3){26} }
\put(-0.8,-0.7){$y$}
\put(16.4,14.3){$d$}
\put(-0.8,0.3){2}
\put(-0.8,1.3){3}
\put(-0.8,2.3){4}
\put(-0.8,3.3){5}
\put(-0.8,4.3){6}
\put(-0.8,5.3){7}
\put(-0.8,6.3){8}
\put(-0.8,7.3){9}
\put(-1.1,8.3){10}
\put(-1.1,9.3){11}
\put(-1.1,10.3){12}
\put(-1.1,11.3){13}
\put(-1.1,12.3){14}
\put(-1.1,13.3){15}
\put(0.4,14.3){1}
\put(1.4,14.3){2}
\put(2.4,14.3){3}
\put(3.4,14.3){4}
\put(4.4,14.3){5}
\put(5.4,14.3){6}
\put(6.4,14.3){7}
\put(7.4,14.3){8}
\put(8.4,14.3){9}
\put(9.1,14.3){10}
\put(10.1,14.3){11}
\put(11.1,14.3){12}
\put(12.1,14.3){13}
\put(13.1,14.3){14}
\put(14.6,14.3){$x$}
\put(0.4,6.3){7}
\put(0.4,7.3){8}
\put(0.4,8.3){9}
\put(0.1,9.3){10}
\put(0.1,10.3){11}
\put(0.1,11.3){12}
\put(0.1,12.3){13}
\put(0.1,13.3){14}
\put(1.1,6.3){13}
\put(1.1,7.3){15}
\put(1.1,8.3){17}
\put(1.1,9.3){19}
\put(1.1,10.3){21}
\put(1.1,11.3){23}
\put(1.1,12.3){25}
\put(1.1,13.3){27}
\put(2.1,6.3){18}
\put(2.1,7.3){21}
\put(2.1,8.3){24}
\put(2.1,9.3){27}
\put(2.1,10.3){30}
\put(2.1,11.3){33}
\put(2.1,12.3){36}
\put(2.1,13.3){39}
\put(3.1,6.3){22}
\put(3.1,7.3){26}
\put(3.1,8.3){30}
\put(3.1,9.3){34}
\put(3.1,10.3){38}
\put(3.1,11.3){42}
\put(3.1,12.3){46}
\put(3.1,13.3){50}
\put(4.1,6.3){25}
\put(4.1,7.3){30}
\put(4.1,8.3){35}
\put(4.1,9.3){40}
\put(4.1,10.3){45}
\put(4.1,11.3){50}
\put(4.1,12.3){55}
\put(4.1,13.3){60}
\put(5.1,6.3){27}
\put(5.1,7.3){33}
\put(5.1,8.3){39}
\put(5.1,9.3){45}
\put(5.1,10.3){51}
\put(5.1,11.3){57}
\put(5.1,12.3){63}
\put(5.1,13.3){69}
\put(6.1,6.3){28}
\put(6.1,7.3){35}
\put(6.1,8.3){42}
\put(6.1,9.3){49}
\put(6.1,10.3){56}
\put(6.1,11.3){63}
\put(6.1,12.3){70}
\put(6.1,13.3){77}
\put(7.1,7.3){36}
\put(7.1,8.3){44}
\put(7.1,9.3){52}
\put(7.1,10.3){60}
\put(7.1,11.3){68}
\put(7.1,12.3){76}
\put(7.1,13.3){84}
\put(8.1,8.3){45}
\put(8.1,9.3){54}
\put(8.1,10.3){63}
\put(8.1,11.3){72}
\put(8.1,12.3){81}
\put(8.1,13.3){90}
\put(9.1,9.3){55}
\put(9.1,10.3){65}
\put(9.1,11.3){75}
\put(9.1,12.3){85}
\put(9.1,13.3){95}
\put(10.1,10.3){66}
\put(10.1,11.3){77}
\put(10.1,12.3){88}
\put(10.1,13.3){99}
\put(11.1,11.3){78}
\put(11.1,12.3){90}
\put(12.1,12.3){91}

\put(11.0,13.3){\small{102}}
\put(12.0,13.3){\small{104}}
\put(13.0,13.3){\small{105} }
   \thinlines
\put(2,2){{\line(1,1){5}}}
\put(1,2){{\line(1,1){6}}}
\put(1,3){{\line(1,1){6}}}
\put(1,4){{\line(1,1){10}}}
\put(1,5){{\line(1,1){9}}}
\put(1,6){{\line(1,1){8}}}
\put(1,7){{\line(1,1){7}}}
\put(1,8){{\line(1,1){6}}}
\put(1,9){{\line(1,1){5}}}
\put(1,10){{\line(1,1){4}}}
\put(1,11){{\line(1,1){1}}}
\end{picture}
   \caption{$I(2,15)$ with non-admissible points crossed out}
     \label{fig:u5}
\end{figure}
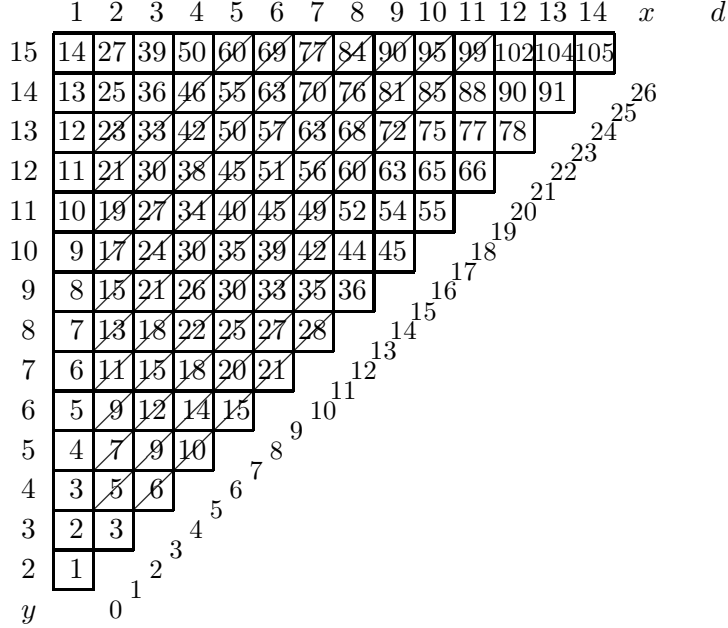

\begin{proof}
(i) The cases $d=0,1,2$ are left to the reader. For $d \ge 3$ the concavity of the cardinality function $c(x,y)$ along the diagonals implies that it is enough to check the admissability of the points $(2,d+1)$, and in
addition $(\frac{d+2}{2},\frac{d+4}{2})$ when $d$ is even, and $(\frac{d+3}{2},\frac{d+5}{2})$ when $d$ is odd. We compare their $c$-values with that of $(1,d+1)$ on the diagonal above. We obtain: 
$$c(2,d+1)=2d-1>d+1=c(1,d+3).$$
Furthermore 
$c(\frac{d+2}{2},\frac{d+4}{2})=\frac{d^2+6d+8}{8}>d+1=c(1,d+3)$ when $d$ is 
even, and 
$c(\frac{d+1}{2},\frac{d+5}{2})=\frac{d^2+8d+7}{8}>d+1=c(1,d+3)$ when $d$ is 
odd.

(ii) Again, by the concavity of the function $c(x,y)$ restricted to 
the diagonal $D_{d}$, 
it suffices to compare a 
few points on each diagonal with the end points on the next.  
Let us first consider the case $d$ even. There is only something to prove if there are at least
$5$ points on $D_d$. The number of points on $D_d$ is 
$\frac{d+2}{2}-(d-m+3)+1=\frac{-d}{2}+m-1$. This is at least $5$ if and only if
$d \le 2m-12$ and $m \ge 9$. So we assume that. 
To show the assertion we must show that for each of the values $c(d-m+5,m-2)$ and $c(\frac{d-2}{2},\frac{d+8}{2})$ (for the points third from the left and third from the right on $D_d$), they cannot be smaller than both of the values $c(d-m+4,m)$ and
$c(\frac{d+2}{2},\frac{d+6}{2})$ for the endpoints of the diagonal $D_{d+1}$.

 For the point $(d-m+5,m-2)$ we obtain two inequalities that reduce to
 $4m-15<3d$, and 
 $$(2m-10-d)(5d-6m+20)<16.$$
 To satisfy this latter inequality, either the factors have different sign or 
 both factors have small positive value.  Furthermore we know that
 $d \le 2m-12$, so $2m-10-d>0$, so only the second factor can be 
 negative.  In this case,
 $4m-15<3d$ and $5d-6m+20<0$, which gives $m\leq 7$, contrary to $m \ge 9$. 
 Thus both factors must have small (even) value.
 If the first factor $2m-10-d$ is $2$, i.e. $d=2m-12$ we must have $5d-6m+20\le 6$,
 which gives $m \le 11$ (and $m \ge 9$. A quick check gives $(4,9)$ as 
 the only  admissible point in this case. 
 A further check of the other cases    $2m-10-d=4$, and $2m-10-d=6$, only possible for $m \le 13$ and $m \le 15$ reveals that the points $(d-m+5,m-2)$ in question are non-admissible in these cases.

The point $(\frac{d-2}{2},\frac{d+8}{2})$ is admissible only if
$$c(\frac{d-2}{2},\frac{d+8}{2})<c(d-m+4,m)\quad {\rm and} \quad 
c(\frac{d-2}{2},\frac{d+8}{2})<c(\frac{d+2}{2},\frac{d+6}{2}).$$
The latter inequality yields $d<12$, and hence $m<d+3=15$, in which 
case the first inequality is satisfied only if  $m=11$ 
and $d=10$, when the point $(4,9)$ is also admissible.

If $d$ is odd, we compare the value of $c$ at the 
 endpoints $(d-m+4,m)$ and $(\frac{d+3}{2},\frac{d+5}{2})$  
on $D_{d+1}$, with its value at $(d-m+5,m-2)$ and 
$(\frac{d-1}{2},\frac{d+7}{2})$, namely the points number
 three and two from the endpoints on $D_{d}$. 
 We only have an issue if we have at least $4$ elements on the diagonal $D_d$, which gives $d \le 2m-11$, and $m \ge 8$.

 The point $(d-m+5,m-2)$ is admissible only if
 $$c(d-m+5,m-2)<c(d-m+4,m)\quad {\rm and}\quad 
 c(d-m+5,m-2)<c(\frac{d+3}{2},\frac{d+5}{2}).$$
  The first inequality reduces to
 $d > \frac{4m-15}{3}$, while the second one yields 
 $$5d^2+(68-16m)d+12m^2-100m+215>0.$$
% hence
%$$(10m-5d-23)(5d-6m+45)<40.$$
For each fixed $m$ this last inequality has a solution only for d outside an interval
$[d_1,d_2]$, where $d_1 < \frac{4m-15}{3}$ and $d_2 > 2m-11,$ at least for $m \ge 11$.
For lower $m$ one checks the statement case by case.
 
  The point $(\frac{d-1}{2},\frac{d+7}{2})$ is admissible only if
$c(\frac{d-1}{2},\frac{d+7}{2}) < c(\frac{d+3}{2},\frac{d+5}{2}), $ which gives 
$d <7$, and $m < d-3$, which is impossible for $m \ge 8$.

(iii) The point $(x,m)$ with $x<m$ is the endpoint with minimal $x$ 
of the diagonal $D_{d}$  
with $d=x+m-3$.   It has cardinality 
$xm - \frac{x(x+1)}{2}$, and for $x < m-3$ it is admissible only if the cardinality of 
the rightmost endpoint of the diagonal $D_{x+m-2}$ is strictly bigger.

If $x+m$ is odd, then the other endpoint of $D_{x+m-2}$ is
$(\frac{x+m-1}{2},\frac{x+m+3}{2})$ and the inequality becomes
$$xm - 
\frac{x(x+1)}{2}<\frac{x+m-1}{2}\frac{x+m+3}{2}-
\frac{1}{2}\frac{x+m-1}{2}\frac{x+m+1}{2}$$
which means
$$5x^2+(8-6m)x+(m^2+4m-5) > 0$$
For each fixed $m$ this last inequality has a solution only for $x$ outside an interval
$[x_1,x_2]$, where $x_1 < \frac{m}{5}+3,$ and $x_2 > m-4$.
The value $x=\frac{m}{5}+2$ satisfies the inequality only if $m \le 13$, and 
for these low values of $m$ we see that the result holds.

If $x+m$ is even, then this endpoint is
$(\frac{x+m}{2},\frac{x+m+2}{2})$ and the condition for admissibilty is
$$xm - \frac{x(x+1)}{2}<\frac{x+m}{2}\frac{x+m+2}{2}-
\frac{1}{2}\frac{x+m}{2}\frac{x+m+2}{2}=\frac{1}{8}(x+m)(x+m+2)$$
which gives:
$$5x^2+(6-6m)x+(m^2+2m) > 0.$$
For each fixed $m$ this last inequality has a solution only for $x$ outside an interval
$[x_1,x_2]$, where $x_1 < \frac{m}{5}+2,$ and $x_2 > m-3$.
The value $x=\frac{m}{5}+1$ satisfies the inequality only if $m \le 13$, and 
for these low values of $m$ we see that the result holds.

(iv) In this case the computation is similar. The point $(x,m-1)$ in
the next to upper-row of $D_d$ has cardinality $c((x,m-1))=x(m-1) - 
\frac{x(x+1)}{2}$ and is
admissible only if it has lower cardinality
than the lower endpoint of the diagonal $D_{d+1}$ above. 
If $x+m$ is odd, then the lower end of the diagonal above is
$(\frac{x+m-1}{2},\frac{x+m+1}{2})$,
and its cardinality is $\frac{(x+m+1)(x+m-1)}{8}$. Now the condition
$$x(m-1) - \frac{x(x+1)}{2} < \frac{(x+m+1)(x+m-1)}{8}$$ 
translates to:
$$(x-m+3)(5x-m-3)+8 > 0.$$
The first factor is negative unless $m-3\leq x\leq m-1$, while 
the second factor is negative when $x<\frac{1}{5}(m+3)$.
If we set $x=m-5$, 
the inequality is satisfied only if $m<8$.
Likewise, if  $x=\frac{m}{5}+1$ the inequality is satisfied only if 
$m<10$.  By the concavity argument the result follows.

If $x+m$ is even, then the lower end of the diagonal above is
$(\frac{x+m-2}{2},\frac{x+m+2}{2})$,
and its cardinality is $\frac{(x+m-2)(x+m+4)}{8}$.
The necessary condition for $(x,m-1)$ being admissible is 
 $$x(m-1) - \frac{x(x+1)}{2} < \frac{(x+m-2)(x+m+4)}{8}.$$
This becomes $$(x-m+4)(5x-m-6)+16 > 0.$$
We insert $m-6$ which is the largest $x$-value smaller than $m-4$
making
$x+m$ even and obtain $m\leq 10$.
Likewise, we insert $x=\frac{m}{5}+2$ and obtain that the inequality then holds for $m \leq 12$.
Hence the statement of the lemma holds for $m \geq 13$.
A special check reveals that it holds for $m=11,12$ also.
\end{proof}

We return to the proof of Proposition \ref{leftright} and assume that 
$S_U$ is a Schubert 
union with spanning dimension $K$, and 
that $S_U$ has the maximal number of points among such unions, i.e. 
$g_{U}$ is maximal 
in the lexicographical order.
Therefore the grid $I_{U}$ contains an admissible point 
$\alpha=(x,y)$ in the 
$d(K)$-diagonal, i.e. $x+y-3=d=d(K)$. 
By Lemma \ref{admiss} it suffices to study the
following eight cases.

(a) $d\leq m-3$

(b)  $m>10$ and $\alpha=(d-m+3,m)$ with and $2\leq d-m+3 \le 
\frac{m}{5}+2$, i.e. $m-1\leq d \le \frac{6m}{5}-1$.

(c) $m>10$ and $\alpha=(d-m+4,m-1)$, with $2\leq d-m+4 \le 
\frac{m}{5}+2$, i.e. $m-2\leq d \le \frac{6m}{5}-2$.

(d) $d$ is even and $\alpha= (\frac{d+2}{2},\frac{d+4}{2})$.

(e) $d$ is odd and $\alpha= (\frac{d+1}{2},\frac{d+5}{2})$.

(f) $d$ is even and $\alpha= (\frac{d}{2},\frac{d+6}{2})$.

(g) $m=11$, and $I_U$ intersects the $d(K)$-diagonal in
$(4,9)$.

(h) $m\leq 10$.

In each case we consider the residual grid $\Delta=I_{U}\setminus 
I_{\alpha}(2,m)$, and find the diagonal $D_{d^{\prime}}$ with largest 
$d^{\prime}$ that $\Delta$ intersects.  By the lexicographical 
ordering 
of $g_{U}$, the value of $d^{\prime}$ is determined in a similar 
fashion as $d=d(K)$ by the cardinality of $\Delta$ and the shape of 
the grid $I(2,m)\setminus I_{\alpha}(2,m)$.
Notice that $S_U$ is a finite union of irreducible components, all of them
Schubert varieties.  Furthermore, the point $\alpha$ corresponds to a 
Schubert 
variety component $S_{\alpha}$ of maximal Krull dimension in $S_U$, and 
any point $\beta\in\Delta\cap D_{d^{\prime}}$ corresponds 
to a Schubert variety  $S_{\beta}$ of maximal Krull-dimension among the 
rest of the irreducible components of $S_U$.
First of all the cardinality of $\Delta$ is $e=K-c(\alpha)$, and by 
definition of $d=d(K)$, the Krull dimension $d^{\prime}$ of $S_{\beta}$ is at most 
$d(K)$.
%Since $g_{U}$ is maximal, $\Delta$ reaches the maximal diagonal $D_{d(K,\alpha)}$ 
%among all diagonals reached by residuals $\Delta^{\prime}=G_{U^{\prime}\setminus G_{\alpha}$, where 
%$U^{\prime}$ is a Schubert union of spanning dimension $K$ that contains the Schubert cycle $S_{\alpha}$.
%Notice that $d(K,\alpha)$ by definition is the maximal Krull dimension among the Schubert cycle components of 
%$U\setminus S_{\alpha}$.

Now we treat the 8 cases separately:
In case (a), when $d\leq m-3$, then $\alpha=(2,3)$ or $(1,d+2)$.  If $\alpha=(2,3)$,
then $d(K)=2$, so $K=3$ and $e=K-c(\alpha)=0$, and $S_U=S_{\alpha}$, which is a
Schubert union of type $S_R$.
Consider the case $\alpha=(1,d+2)$. If $K\geq d+2$, then $d(K)>d+1=c(\alpha)$, 
contrary to the assumption, so $K=c(\alpha)=d+1$.
In particular $e=K-c(\alpha)=0$ and $S_U=S_{\alpha}$, which is of type $S_L$.

In case (b),  the grid $I(2,m)\setminus I_{\alpha}(2,m) = \{(x,y) \in I(2,m)| 
d-m+4\leq x<y\leq m\}$.  
Since $d(K)=d$, the cardinality $e=K-c(\alpha)$ of $\Delta$ is less than the 
cardinality of the leftmost column of $I(2,m)\setminus I_{\alpha}(2,m)$, 
i.e. at most $m-(d-m+5)=2m-d-5$.  
We use an argument of (a), almost identical to that in the $\alpha=(1,d+2)$ case of (a),
to conclude that $\Delta=\{(d-m+4,y)| d-m+4<y\leq e\}$.  Notice 
furthermore that $S_U$ clearly is of type $S_L$.

In case (c), with $\alpha=(d-m+4,m-1)$. Since $d(K)=d$, we first see 
that the cardinality 
of $\Delta$ is less than the cardinality of the upper row of 
$I_{\alpha}(2,m)$, i.e.
$e=K-c(\alpha)) < d-m+4$.  Compare now the row of points $R=\{(x,m)| 
1\leq 
x<d-m+4\}$ with the column $C=\{(d-m+5,y)| d-m+5<y<m\}$, both in 
$I(2,m)\setminus I_{\alpha}(2,m)$.
Notice that both have cardinality at least $e$, so that for $\Delta$ 
to 
reach the maximal diagonal $D_{d^{\prime}}$, it must be contained in 
one of these.
 The row $R$ starts on the diagonal 
$D_{m-2}$, while the columns $C$ starts on the diagonal 
$D_{2d-2m+8}$.  When $m>10$ and $d \le \frac{6m}{5}-2$, the highest 
of these diagonals is $D_{m-2}$, since then $2d-2m+8\leq \frac{2m}{5}+4 <m-2$, 
so in that case $\Delta$ must be completely contained in the row $R$.

To see whether $S_U$ is of type $S_L$, there are essentially two different situations:
$e=d-m+3$ (maximum possible), and $e \le d-m+2$. If $e=d-m+3$, we are
already in case b), since $(d-m+3,m)$ and $(d-m+4,m-1)$ are on the 
same diagonal, and the point $(d-m+3,m)$ is covered by case (b).

Assume $e \le d-m+2$. Then $S_U$ is not of type $S_L$, but we revise $I_U$, and collectively remove the $d-m+3-e$ top 
squares
of the right column of $I_{\alpha}(2,m)$, and reinstall them horizontally
as points $(e+1,m), (e+2,m)...,(d-m+3,m)$. This amounts to moving 
squares along  $d-m+3-e$ diagonals, and does not alter the number
of ${\Fq}$-rational points of the Schubert unions represented
by the two grids. But after moving, we have the grid of type $S_L$, as in case (b) and we are 
done. 
(We have $d-m+3$ columns to the left filled up completely,
and $K-c(d-m+3,m)$ squares in column nr $d-m+4$).

In case (d), with $\alpha= (\frac{d+2}{2},\frac{d+4}{2})$  
the lowest row of $I(2,m)\setminus I_{\alpha}(2,m)$, if any, has larger length 
than $d^{\prime}$, since $d(K)>d^{\prime}$.  Therefore $\Delta$ is 
completely contained in this lowest row.  Furthermore, $S_U$ is clearly 
of type $S_R$.

In case (e) we see that $K-c(\alpha)=0$, since if $\Delta$ is 
non-empty, then it could lie to the right of $\alpha$ in the top row, 
and $d(K)$ would have been larger. On the other hand $S_U=S_{\alpha}$ 
is 
clearly of type $S_R$.

In case (f), we see that $K-C(d(K))=0$ or $1$, since if it was
at least $2$, $\Delta$ could lie to the right of $\alpha$, and
then $d(K)$ would have been larger.
If $\Delta$ is empty, there is nothing to do, and if $\Delta$ consists of one point, 
the choice with the largest $g_U$ is  $\Delta=\{(\frac{d+2}{2},\frac{d+4}{2})\}$(lies to the right of and below $\alpha$.
In both cases $S_U$ is of type $S_R$.

Starting with (g), we have $m=11$, and $\alpha=(4,9)$. Since $\alpha$
lies on the diagonal $D_{10}$ and
$c(4,9)+1=27=C(11)$, we see that this is only an issue when 
$K=c(4,9)=26$, i.e. $\Delta$ is empty. But a quick calculation reveals 
that
$S_{(4,9)}$ is not optimal at all for $K=26$ with respect to $g_U$. An optimal choice in this
case is  $S_U=S_{(2,11)} \cup S_{(3,10)}$, 
%with $M_U=\{10,9,7\}$, 
which is of type $S_L$. Hence the assertion holds in all cases.

The cases $m\leq 10$ are entirely similar. The check that no case 
occurs that is not of one of the kinds above is left to the reader.
\end{proof}

It seems obvious that for each given $K$ and $m$ the one(s) among $S_L$ and
$S_R$ which is maximal with respect to the lexicographical order on $g_U$,
has a maximal number of ${\Fq}$-rational points, not only for large enough $q$,
but for all fixed prime powers $q \ge 2$. 
All examples indicate that, but we have not been able to prove it.

% for a given $S_U$
%with spanning dimension $K$ one can obtain the grid corresponding to the $g_U$-maximal,
%say $S_R$, from $I_U$ by moving individual squares, such that no square is moved to
%a diagonal with a lower $d$-value than that of the diagonal it is moved from.

\begin{rem} \label{decomp}
{\rm We see that Proposition \ref{leftright} implies the result hoped for in Remark 29
of \cite{GPP}, at least for large enough $q$. This result says that for the Grassmann code $G(2,m)$, with $m >4$ we have
$$d_r=(q^\delta +q^{\delta-1}+\cdots +q^{\delta -m+2}) +
(q^{\delta-2}+q^{\delta-3}+\cdots +q^{\delta-r+m})$$
and 
$$ d_{k-r}=n-(1+q+\cdots +q^{m-2}) -
(q^{2}+q^{3}+\cdots +q^{r-m+2}).$$
for $m <r\le 2m-5.$ 

This follows from the fact that for $ m < n- K \le 2m-5$ the $S_R$ are maximal with respect to the lexicographical order on the $g_U$, among those Schubert unions with spanning dimension $K$, and for $ m < K \le 2m-5$ the $S_L$ are maximal with respect to the lexicographical order on the $g_U$, among those Schubert unions with spanning dimension $K$. One then simply counts the number of ${\Fq}$-rational points on these
particular unions $S_L$ and $S_R$.

}
\end{rem}


\begin{thebibliography}{[E-L-M-S]}

%\bibitem[GL]{GL}  Ghorpade, S, Lachaud, G. \textit{ Higher weights
%of Grassmann Codes}, in \textit{Coding Theory, Cryptography, and 
%Related
%Areas} (Guanajuoto, 1998), Springer Verlag, Berlin/Heidelberg, 
% 122-31 (2000).

%\bibitem[Fu]{Fu}  W. Fulton,  \textit{Young Tableaux}, 
% Student Texts, 35, London Math. Soc. (1991).

\bibitem[GPP]{GPP}  S. R. Ghorpade, A. R. Patil, H.K. Pillai \textit{ Decomposable subspaces, linear sections of Grassmann varieties, and higher weights of Grassmann codes}, Finite Fields Appl., 15, 54-68 (2009).

%\bibitem[GT]{GT} Ghorpade, S., Tsfasman, M. \textit{Schubert 
%Varieties,
%Linear Codes, and Enumerative Combinatorics}, Preprint (2003).

%\bibitem[HC]{HC}  H. Chen,  \textit{ On the minimum distance of 
%Schubert 
%codes}, IEEE Trans. of Inform. Theory, 46, 1535-38 (2000).


\bibitem[HJR]{HJR}  J. P. Hansen, T. Johnsen, K. Ranestad, 
  \textit{Schubert unions in Grassmann varieties}, Finite Fields Appl., 13, 738-750 (2007).
  
\bibitem[HJR2]{HJR2}  J. P. Hansen, T. Johnsen, K. Ranestad, 
  \textit{Schubert unions in Grassmann varieties}, arXiv: math.AG/0503121, (2005).


%\bibitem[N]{N} Nogin, D.Yu., \textit{Codes associated to 
%Grassmannians}, 
%in \textit{Arithmetic Geometry and Coding theory} (Luminy 1993),
%R. Pellikaan, M. Perret, S.G. Vladut, Eds. Walter de Gruyter,
%Berlin/New York,  145-54 (1996).

\bibitem[P]{P}  A. R. Patil,  \textit{Weight hierarchy and generalized spectrum of linear codes associated to Grassmann varieties}, Ph.D. thesis, Indian Institute of Technology
Bomabay (2008).


\end{thebibliography}
\end{document}